\documentclass[english,a4paper,11pt]{amsart}
\usepackage{srcltx}
\usepackage{amsmath}
\usepackage{amsfonts}
\usepackage{mathrsfs}
\usepackage[all]{xy}
\usepackage{amsthm}
\usepackage{latexsym}
\usepackage{amssymb}
\usepackage{graphicx}
\usepackage{verbatim}
\usepackage{fancyhdr}
\usepackage{hyperref}

\newtheorem{lemma}{Lemma}[section]
\newtheorem{theorem}[lemma]{Theorem}
\newtheorem{proposition}[lemma]{Proposition}
\newtheorem{corollary}[lemma]{Corollary}

\theoremstyle{definition}
\newtheorem{definition}[lemma]{Definition}
\theoremstyle{remark}
\newtheorem{remark}[lemma]{Remark}

\newtheorem{example}[lemma]{Example}

\newcommand{\C}[0]{\mathbb{C}}

\newcommand{\p}[0]{\mathbb{P}}

\def\cocoa{{\hbox{\rm C\kern-.13em o\kern-.07em C\kern-.13em o\kern-.15em A}}}

\begin{document}

\title{Apolarity, Hessian and Macaulay polynomials}
\author{Lorenzo Di Biagio}
\address{Dipartimento di Matematica, Universit\`a degli Studi ``Roma Tre'' - Largo San Leonardo Murialdo 1, 00146, Roma, Italy}
\email[1]{dibiagio@mat.uniroma3.it}
\email[2]{lorenzo.dibiagio@gmail.com}
\author{Elisa Postinghel}
\address{Centre of Mathematics for Applications, University of Oslo - P.O. Box 1053 Blindern, N0-0316 Oslo, Norway}
\email[1]{elisa.postinghel@cma.uio.no}
\email[2]{elisa.postinghel@gmail.com}
\maketitle

\begin{abstract} 
A result by Macaulay states that an Artinian graded Gorenstein ring $R$ of socle dimension one and socle degree $\delta$ can be realized as the apolar ring $\mathbb{C}[\frac{\partial}{\partial{x_0}}, \ldots, \frac{\partial}{\partial{x_n}}]/{g^\perp}$ of a homogeneous polynomial $g$ of degree $\delta$ in $x_0,\dots,x_n$. If $R$ is the Jacobian ring of a smooth hypersurface $f(x_0,\dots,x_n)=0$ 
then $\delta$ is equal to the degree of the Hessian polynomial of $f$. In this paper we investigate the relationship between $g$ and the Hessian polynomial of $f$ and we provide a complete description for $n=1$ and $\deg(f)\leq4$ and for $n=2$ and $\deg(f)\leq3$.
\end{abstract}

\section{Introduction: the problem}

Let $f \in R=\mathbb{C}[x_0, \ldots, x_n]$ be a homogeneous polynomial. The ideal $J(f)$ generated by the partial derivatives of $f$ is called the Jacobian ideal, or the gradient ideal, of $f$. In the smooth case this ideal contains a power of the irrelevant ideal, so it has maximum depth in the coordinates ring and it is generated by a regular sequence. The associated ring $R(f)=R/J(f)$, the so-called Jacobian ring of $f$, is an Artinian Gorenstein graded ring. Formal definitions can be found in Section \ref{Hessiano e Macaulay}. Both the Jacobian ideal and its counterpart, the Jacobian ring, have been largely studied and it is now clear, form the work of P.\ Griffiths, that they reflect many geometric properties of the variety. For example if $f$ and $f'$ define smooth hypersurfaces then $f$ and $f'$ are projectively equivalent if and only if $R(f)$ is isomorphic to $R(f')$ thus allowing us to recover $V(f)$ from its Jacobian ring. Moreover partial information about $R(f)$ is equivalent to information about the Hodge groups that appear in the Hodge decomposition of $H^{n-1}(V(f),\mathbb{C})$.  See \cite[\S 2]{C} for a nice account on this stuff and further references. Neverthless so far the Jacobian ring has not been completely understood. 

Apolarity allows us to associate an Artinian Gorenstein graded  ring to a form. This nice property is described in a classical theorem due to Macaulay (Theorem \ref{Macaulay}): there exists a homogeneous polynomial (the Macaulay polynomial) $g$ such that $J(f)$ is equal to $g^\perp$, where $g^\perp \subset T=\mathbb{C}[\frac{\partial}{\partial{x_0}}, \ldots, \frac{\partial}{\partial{x_n}}]$ (upon identifying $x_i$ with $\frac{\partial}{\partial x_i}$).
Apolarity is a very good tool to study varieties of sum of powers (see for example \cite{IR2, IR1, RS}), which are the objects of a very deep and challenging research area in Algebraic Geometry.

It is not so immediate to compute by hand the Macaulay polynomial associated to a given Artinian Gorenstein graded ring, but in the case of the Jacobian ring it seems natural to look at the Hessian polynomial $\text{Hess}(f)$ of $f$ since it has the right degree. It immediately turns out that if $f$ is a Fermat polynomial, then $\text{Hess}(f)$ and the Macaulay polynomial associated to $R(f)$ coincide, up to scalars (see Example \ref{fermat}).
Therefore we ask ourselves if $\text{Hess}(f)$ is always the Macaulay polynomial (up to scalar multiplication) associated to $f$.
A first naive conjecture is the following:
$J(f) = {\text{Hess}(f)^\perp}$, for every smooth homogeneous polynomial $f\in R$.

We will  see in Section \ref{Hessiano e Macaulay} that the question is not actually meaningful, anyway the answer is `no' in general, but `yes' in certain cases.

In Section \ref{binary} we will study the question for binary forms, giving a complete answer for forms of degree $3$ and $4$.

In Section \ref{cubiche piane}
we will completely answer the question in the plane cubics case. 

In Section \ref{quartiche piane} we will present a first attempt to study plane quartics, giving some specific examples.

In Section \ref{cocoa} we will see how to use the computer algebra system \cocoa \ to attack this problem.

\section{Preliminaries}

We work over the complex numbers $\C$, but all the results hold for any algebraically closed field of characteristic zero. Let $S:=\C[x_0,\ldots,x_n]$ be the polynomial ring in $n+1$ variables, let $T:=\C[\partial_0, \ldots, \partial_n]$ be the $\C$-algebra generated by the partial derivatives $\partial_i$, where $\partial_i:=\frac{\partial}{\partial x_i}$. $S$ and $T$ are naturally graded rings and we denote by $S_d$ and $T_d$ their degree $d$ part, which is of course a $\C$-vector space of dimension $\binom{n+d}{d}$.

By the natural differentiation action of $T$ on $S$ we can view $S$ as a $T$-module. 
Analogously we can think of $S$ as the algebra of partial derivatives on $T$, hence we can also view $T$ as an $S$-module. These two actions define a perfect pairing between homogeneous forms of degree $j$ (cf.\ \cite[Prop.\ 2.3]{G} or \cite{IK}):
\begin{equation} \label{pairing}
S_j \times T_j \rightarrow \C.
\end{equation}


Given $g \in S$ and $f \in T$ we will say that $f$ is \emph{apolar} to $g$ if $f \cdot g = 0$. 

\begin{remark} \label{commutativity}
If $f,g$ are homogeneous of the same degree then $g \cdot f= f \cdot g$. 
\end{remark}

Regarding $S$ as a (left) $T$-module, let $I \subseteq T$ be an ideal and let $M, P \subseteq S$ be $T$-submodules of $S$.  Recall that  $(P:_IM):=\{i \in I | iM \subseteq P \}$ is an ideal of $T$ contained in $I$. If $M$ is principal, $M=Tg$, then we will write $(P:_Ig)$ instead of $(P:_IM)$. Analogously, if $M \subseteq P$, then recall that $(M:_PI):=\{p \in P | ip \in M, \forall i \in I \}$ is a $T$-submodule of $P$. 

In particular for any polynomial $g \in S \setminus \{0\}$ we will denote by $g^\perp$ the ideal of $T$ of forms apolar to $g$, i.e. $g^\perp:=\textrm{Ann}(g)=\{f \in T | f \cdot g = 0\}=(0:_Tg)$. Let $T_g:={T}/{g^\perp}$; since $\sqrt{g^\perp}=(\partial_0, \ldots, \partial_n)$ then $T_g$ is an Artinian ($0$-dimensional) local ring. 

Recall that a zero-dimensional local ring $A$ is Gorenstein if and only if its \emph{socle} (i.e. the annihilator of the unique maximal ideal) is simple (cf.\ \cite[Prop.\ 21.5]{Ei}). If moreover the ring is graded, we will call \emph{socle degree} the maximum integer $j$ such that $A_j \not = 0$.
Recall the following theorem (see \cite[theorem 21.6]{Ei}, \cite[\S 60ff]{Ma}, \cite[lemma 2.12]{IK} or the lecture notes by Geramita, esp.\ lecture 8 in \cite{G}):

\begin{theorem}[Macaulay] \label{Macaulay}
With notation as above, there is a one-to-one inclusion reversing correspondence between finitely generated nonzero $T$-submodules $M \subseteq S$ and ideals $I \subseteq T$ such that $I \subseteq (\partial_0, \ldots, \partial_n)$ and $T/I$ is a local Artinian ring, given by
$$M \mapsto (0:_TM), \text{ the annihilator of }M\text{ in }T;$$
$$I \mapsto (0:_SI), \text{ the submodule of }S\text{ annihilated by }I.$$
 
In particular, ideals $I$ as above such that  $T/I$ is local Artinian Gorenstein correspond to principal submodules $Tg$ for some element $g \in S \setminus\{0\}$ (i.e.\ $I=g^\perp$).
\end{theorem}
\begin{corollary}[Macaulay] \label{Macaulayomogeneo}
The homogeneous ideals $I$ as in Theorem \ref{Macaulay} such that $T/I$ is graded local Artinian Gorenstein and of socle degree $j$ correspond to principal submodules $Tg$, where $g$  is a homogeneous polynomial in $S_j$.
\end{corollary}



 
 \begin{definition}
We will call \emph{Macaulay polynomial} associated to $T/I$ the polynomial $g$ associated to $T/I$ (up to scalar multiplication) as in Corollary \ref{Macaulayomogeneo}.  
 \end{definition}
 
 \begin{remark}
Some authors refer to the \emph{Macaulay polynomial} as the \emph{dual socle generator} of $T/I$.
 \end{remark}
\begin{remark}
If $g$ is a homogeneous polynomial of degree $j$, then the Hilbert function $h(T_g)$ is symmetric with respect to $j/2$ (cf. \cite[p.\ 9]{IK}). Hence by Corollary \ref{Macaulayomogeneo}, if $I \subseteq (\partial_0, \ldots, \partial_n)$ is a homogeneous ideal of $T$ and $T/I$ is (graded) local Artinian Gorenstein and of socle degree $j$, then its Hilbert function $h(T/I)$ is symmetric with respect to $j/2$.
\end{remark}

Given $I \subseteq (\partial_0,\ldots ,\partial_n)$ homogeneous ideal such that $T/I$ is local Artinian Gorenstein and of socle degree $j$, one can wonder if there is a simple way to determine the associated Macaulay polynomial $g \in S_j$. In fact the following holds:
\begin{remark} \label{osservazione}
The Macaulay polynomial associated to $T/I$ is any nonzero element of $(0:_{S_j}I_j)$.
\end{remark} 

\begin{proof}
By Corollary \ref{Macaulayomogeneo}, we know that $(0:_SI)$ is a principal $T$-submodule generated by $g$, hence $$(0:_SI)=Tg=T_0g \oplus \cdots  \oplus T_jg=(0:_{S_j}I) \oplus \cdots \oplus (0:_{S_0}I),$$ and $T_0g=\C g$, the $\C$-vector space of dimension $1$ generated by $g$. Therefore any nonzero element of $T_0g=(0:_{S_j}I)=\{s \in S_j | is=0 \ \forall i \in I \}$ can be chosen as $g$.  Moreover since $T/I$ has socle degree $j$ and socle dimension $1$, then $I_j$ is a $\C$-vector subspace of $T_j$ of codimension $1$, hence, by (\ref{pairing}), $(0:_{S_j}I_j)$ has dimension $1$. Since $(0:_{S_j}I_j) \supseteq (0:_{S_j}I) $ and they have the same dimension we have that $(0:_{S_j}I_j) = (0:_{S_j}I)$.
\end{proof}


\section{Jacobian ring and Hessian polynomial}\label{Hessiano e Macaulay}
Set $R:=\C[x_0, \ldots, x_n]$. Let $R_d$ be its homogeneous degree $d$ part. Clearly $R=S$, but when we write $R$ instead of $S$ we stress the fact that we view the polynomial ring as the  $\C$-algebra of partial derivatives, by the action $f(x_0, \ldots, x_n) \cdot g(x_0, \ldots, x_n)=f(\partial_{0}, \ldots, \partial_n)(g(x_0, \ldots, x_n))$, hence implicitly identifying $R$ with $T$. Let $f \in R_d$ defining a smooth hypersurface $V(f) \subset \mathbb{P}^n=\bold{Proj}(R)$. Let $J(f)$ be the \emph{Jacobian ideal} of $f$, i.e., the homogeneous ideal in $R$ generated by the partial derivatives $\frac{\partial f}{\partial x_0}, \ldots, \frac{\partial f}{\partial x_n}$. Since $f$ is smooth, then $\sqrt{J(f)}=(x_0, \ldots, x_n)$, hence the \emph{Jacobian ring} $R(f):=R/J(f)$ is graded, Artinian and local. It is well-known that in this case $R(f)$ is a complete intersection hence it is also Gorenstein and of socle degree $(n+1)(d-2)$ (see \cite{Ei}).

By the identification of $R$ with $T$, we see that given $f \in R_d$, $f \not = 0$, the Jacobian ring $R(f)$ satisfies the hypotheses of Corollary \ref{Macaulayomogeneo}, i.e., we can associate to $R(f)$ its Macaulay polynomial, namely the homogeneous polynomial $g \in S$ of degree $(n+1)(d-2)$ such that, under the natural identification of $R$ with $T$, $J(f)=g^\perp$. We will call this $g$ the Macaulay polynomial associated to $f$, meaning that $g$ is the Macaulay polynomial associated to $R(f)$. We will denote it by $g=\text{Mac}(f)$.

Given $f(x_0, \ldots, x_n)$ as before, we can also consider the Hessian polynomial of $f$, $\text{Hess}(f) \in S$, that is the determinant of the Hessian matrix, i.e., the matrix of partial, second order derivatives of $f(x_0, \ldots, x_n)$. We will always consider $\text{Hess}(f)$ up to scalar multiplication. Since $f$ is nonsingular, by \cite[\S 2.2]{CRS} we have that $\text{Hess}(f) \not = 0$. Moreover $\text{Hess}(f)$ is homogeneous and $\deg(\text{Hess}(f))=(n+1)(d-2)$. 

As we mentioned in the introduction, we first ask ourselves whether $\text{Hess}(f)$ is the Macaulay polynomial (up to scalar multiplication) associated to $f$. 
Notice that checking if $\text{Hess}(f)$ is the socle generator of the Jacobian ring associated to $f$ is not a difficult task: by Remark \ref{osservazione} and Remark \ref{commutativity} it is necessary and sufficient to check that $\text{Hess}(f)$ kills all the forms in $J(f)_{(n+1)(d-2)}$.

\begin{example} \label{fermat}
Let $f = x_0^d+\cdots + x_n^d \in R$ be the Fermat polynomial of degree $d$ in $n+1$ variables. Then $J(f)=(x_0^{d-1}, \ldots, x_n^{d-1})$ and $\text{Hess}(f)$ is a monomial, $\text{Hess}(f)=(d(d-1))^{n+1} x_0^{d-2} \cdot \cdots \cdot x_n^{d-2}$. In this case $\text{Hess}(f)$ is the Macaulay polynomial associated to $f$. In fact for any monomial $p \in S$ of degree $n(d-2)-1$, $\text{Hess}(f) \cdot x_i^{d-1} p  = 0 \Leftrightarrow \forall c \in \C\setminus \{0\}, \text{Hess}(f)\neq c x_i^{d-1}p$, and this inequality clearly always holds, for any $c$, $p$ and $i$. Thus  $R(f)= T/{\text{Hess}(f)^\perp}$.
\end{example}

\begin{example} \label{esempio}
Let $n=2$, $f=(x_0+x_1)^3+x_1^3+x_2^3.$ Then $J(f)=(3(x_0+x_1)^2, 3(x_0+x_1)^2+3x_1^2, 3x_2^2)$ and $\text{Hess}(f)=216(x_0x_1x_2 + x_1^2x_2)$. In this case $\text{Hess}(f)$ is not the Macaulay polynomial associated to $f$. In fact, for example, $3(x_0+x_1)^2x_2 \in J(f)$, but $\text{Hess}(f) \cdot 3(x_0+x_1)^2x_2 = 2592 \not = 0$. The Macaulay polynomial associated to $f$ is (up to scalars) $x_0^2x_2 - x_0x_1x_2$.
\end{example}

By Example \ref{fermat} and Example \ref{esempio},  since the second one is simply obtained by performing a change of variables in the first one, it should be clear that the question whether $\text{Hess}(f)$ is equal (up to scalars) to $\text{Mac}(f)$ is not the right one.

Identifying $R$ with $T$, as before, and taking $f \in R_d=T_d$, $f \not = 0$, first of all we should understand how $\text{Hess}(f)$ and $\text{Mac}(f)$ behave under a linear change of variables, i.e., under the action of $SL_{n+1}(\C)$ on $R_1$. Let $\overline{x}=(x_0, \ldots, x_n)$. Let $A \in SL_{n+1}(\C)$. 
The Hessian polynomial $\text{Hess}(f)$ is covariant under the change of variables given by $A$, that is: $\text{Hess}(f(A\overline{x}))=\text{Hess}(f(\overline{x}))(A \overline{x})$ (see \cite[\S 2.1]{CRS}). For the Macaulay polynomial the following lemma holds. It is well-known to experts but we include a proof for lack of references: 

\begin{lemma}
$\text{Mac}(f)$ is contravariant under the change of variables given by $A$, that is: $\text{Mac}(f(A\overline{x}))=\text{Mac}(f(\overline{x}))(^tA^{-1}\overline{x})$.
\end{lemma}

\begin{proof}
The two bases $\overline{x}=(x_0, \ldots, x_n)$ of $S_1$ and $R_1$ are dual to each other under the derivation actions (we are identifying $R$ and $T$, as usual). Therefore the new basis $A \overline{x}$ of $R_1$ is dual to the basis $^tA^{-1} \overline{x}$ of $S_1$, hence for any polynomials $p \in S_d$ and $q \in R_d$ we have $p(\overline{x}) \cdot q(\overline{x})=p(^tA^{-1} \overline{x}) \cdot q(A \overline{x})$.

Moreover $J(q(A \overline{x}))$ is equal to $J(q)(A \overline{x})$, in fact $\nabla (q(A \overline{x}))=(\nabla q)(A \overline{x}) \cdot \nabla(A\overline{x})=(\nabla q)(A \overline{x}) \cdot A.$ Since $A$ is invertible, the ideal generated by the entries in $\nabla (q(A \overline{x}))$ is the same as the ideal generated by the entries in $(\nabla q)(A \overline{x}) \cdot A$.
The thesis now follows. 
\end{proof}

Therefore, the new question we are interested in is ``when are $\text{Mac}(f)$ and $\text{Hess}(f)$ projectively equivalent?''.

\begin{example}
Going back to Example \ref{esempio}, as we have noted before, the polynomial $f$ is obtained from the Fermat cubic polynomial $h:=x_0^3+x_1^3+x_2^3$ by the linear change of variables given by the matrix
$$ A:=\left( \begin{array}{ccc}
1 & 1 &0 \\
0 & 1 & 0\\
0 & 0 & 1\\
\end{array} \right). $$
Hence $\text{Hess}(f(\overline{x}))=\text{Hess}(h)(A \overline{x})$, while $\text{Mac}(f(\overline{x}))=\text{Mac}(h)(^tA^{-1}\overline{x})$. Therefore $\text{Mac}(f)$ and $\text{Hess}(f)$ are projectively equivalent by the linear change of variables given by $^tAA$, as it can be easily verified.
\end{example}

\section{Binary forms of degree $3$ and $4$} \label{binary}
In this section we will deal with homogeneous polynomials in two variables. For the sake of simplicity in this and following paragraphs we will use the variables $x,y,z, \ldots$ instead of $x_0,x_1,x_2, \ldots $.

Let $f(x,y)$ be a nonsingular homogeneous polynomial of degree $d$; $f(x,y)=g_1(x,y)\cdot \cdots \cdot g_d(x,y)$, where $g_i(x,y)=a_ix+b_iy$ are linear forms that are distinct up to constants. 

If $d=3$, since any set of three distinct points in $\mathbb{P}^1$ is projectively equivalent to any other set of this type, then every $f$ as before is projectively equivalent to $x^3+y^3$, the Fermat polynomial. Hence $\text{Hess}(f)$ and $\text{Mac}(f)$ are projectively equivalent for any $f$.

If $d=4$, then any $f$, by the same token as before, is projectively equivalent, up to constants, to $f_a(x,y)=xy(x-y)(x+ay)$, where $a \in \C, a \not = 0,-1$. Set $H_a:=\text{Hess}(f_a)$ and $M_a:= \text{Mac}(f_a)$. If $a=1$ then $f_1$ is projectively equivalent to the Fermat polynomial $x^4+y^4$.  In this case we have $$H_1=-9(x^4+2x^2y^2+y^4), \ M_1=x^4+2x^2y^2+y^4$$ and clearly these two (singular) polynomials are equal (up to constants).

From now on we can suppose that $a \not = -2,-\frac{1}{2}, 1$, since the case $a=1$ 
has just been analyzed, and when $a=-2, - \frac{1}{2}$ we have that $f_a$ is projectively 
equivalent to $x^4+y^4$. The Jacobian ideal of $f_a$ is 
$$J(f_a)=\left(3 x^2 y - 2 x y^2 + 2 a x y^2 - ay^3, x^3 - 2 x^2 y + 2a x^2 y - 3 ax y^2 \right).$$
The Hessian polynomial of $f_a$ is 
$$H_a =-9x^4-12(a-1)x^3y-6(2a^2-a+2)x^2 y^2+12 a(a-1)x y^3-9a^2 y^4. $$ 
It can be easily seen that the polynomial 
$$M_a =(a^2+a+1)x^4-2(a-1)x^3y+6x^2 y^2+2\frac{a-1}{a}x y^3+\frac{a^2+a+1}{a^2} y^4$$ 
is the Macaulay polynomial of $f_a$: indeed, it is apolar to the degree $4$ part of $J(f_a)$ that is generated by the forms
$x(3 x^2 y - 2 x y^2 + 2 a x y^2 - ay^3)$,  $x(x^3 - 2 x^2 y + 2a x^2 y - 3 ax y^2)$, $y(3 x^2 y - 2 x y^2 + 2 a x y^2 - ay^3) $ and $y(x^3 - 2 x^2 y + 2a x^2 y - 3 ax y^2)$
 (cf. Remark \ref{osservazione}).

Recall that if $\lambda$ is the cross-ratio of four distinct ordered points $V(f)$ in $\mathbb{P}^1$, then the associated $j$-invariant is defined as $$j(f):=2^8\frac{(\lambda^2-\lambda+1)^3}{\lambda^2(\lambda-1)^2}$$ and it does not depend on the order of the points.
The values of the $j$-invariant correspond to the projective equivalence classes of binary quartic forms.

\begin{proposition}
There are only three projective equivalence classes of smooth binary quartic forms such that the Macaulay and the Hessian polynomial are projectively equivalent. They correspond to the values $0$, $1728$ and $6912$  of the $j$-invariant.
\end{proposition}
\begin{proof}
Recalling that the class of the Fermat quartic has $j$-invariant equal to $1728$, by the preceding arguments we can assume that  $a \neq0,1,-1,-2,-\frac{1}{2}$.

Since $H_a$ and $M_a$ are two homogeneous polynomials of degree $4$, they are projectively equivalent if and only if they have the same $j$-invariant. 
By computing the solutions $V(H_a)$ (and $V(M_a)$) of $H_a=0$ (and $M_a=0$, respectively) and the cross-ratios, we get  \begin{align*} j(H_a)&= 
2^8 \frac{(1+a+a^2)^6}{a^2(a+1)^2(a-1)^2(a+1/2)^2(a+2)^2}, \\ j(M_a)&= 
2^8\frac{27(1+a+a^2)^3}{(a-1)^2(a+1/2)^2(a+2)^2}. \end{align*} By our assumptions on $a$ these numbers are well-defined.
The equation $j(H_a)=j(M_a)$, that is equivalent to the following 
$$ (1+a+a^2)^3 \cdot [(1+a+a^2)^3-27 a^2(a+1)^2]=0, $$
gives us the twelve values of $a  \in \C\setminus\{0,1,-1,-2,-\frac{1}{2}\}$ (counted with multiplicity)  
such that the corresponding binary forms have Hessian and Macaulay polynomials projectively equivalent. 
In particular for six of them, namely the solutions given by the first factor of the equation, we have that $j(f_a)$ is equal to $0$, that $H_a$ and $M_a$ are nonsingular and that $j(H_a)=j(M_a)=0$; for the other six, the ones coming from the second factor of the equation, we have that $j(f_a)=6912$, that $H_a$ and $M_a$ are nonsingular and that $j(H_a)=j(M_a)=2304$. 

In all remaining infinitely many cases $H_a$ and $M_a$ are not projectively equivalent.
\end{proof}

\begin{remark}
Given a nonsingular binary quartic form $f$, notice that  $f$, $\textrm{Hess}(f)$ and $\textrm{Mac}(f)$ are all mutually projectively equivalent if and only if $j(f)=0$.
\end{remark}

\section{Plane cubics}\label{cubiche piane}
In this section we will investigate the relations between the Hessian and the Macaulay polynomial associated to any smooth cubic in $\mathbb{P}^2$. It is the first interesting case to analyze in the plane, in fact in the conic case the two polynomials turn out to be trivially equal (up to scalar multiplication). 

Recall that the \emph{Hasse pencil} is the one-parameter family of curves defined by 
\begin{eqnarray}\label{fascio di Hasse}
f_a(x,y,z)=x^3+y^3+z^3-3axyz, \ a \in \mathbb{C}.
\end{eqnarray}

Any smooth cubic in $\p^2$ is projectively equivalent to a cubic in the Hasse pencil for a certain $ a\in\C$, $a^3\neq 1$. Hence we can reduce the problem of studying planar cubics to the analysis of the cubics in (\ref{fascio di Hasse}). 
If $a=0$ we simply get the Fermat cubic and, as already seen in Example \ref{fermat}, $H_0=M_0=xyz$. From now on we will suppose $a \not = 0$.
The Jacobian ideal of $f_a(x,y,z)$ is 
$$
J(f_a)=\left( x^2-ayz,y^2-axz,z^2-axy\right).
$$
The Hessian polynomial is
$$
H_a:=\textrm{Hess}(f_a)=\left|
\begin{array}{ccc}2x&-az&-ay\\
-az&2y&-ax\\
-ay&-ax&2z
\end{array}\right|
=(8-2a^3)xyz-2a^2 (x^3+y^3+z^3).
$$
Notice that $H_a$ is equal (up to scalar multiplication) to $f_b(x,y,z)$, with
\begin{equation} \label{ahessiano}
b={(4-a^3)}/{3a^2}.
\end{equation}
By Remark \ref{osservazione} we look for the Macaulay polynomial $M_a$ of $f_a$ among all cubic forms apolar to $J(f_a)_3$. We notice that $M_a$ has to be symmetric with respect to $x,y,z$ 
since $J(f_a)_3$ is invariant under permutation of the variables. Therefore, with a simple computation, it is straightforward to see that $M_a$ is equal (up to scalar multiplication) to $f_c$, with
\begin{equation} \label{amacaulay}
c=-{2}/{a}.
\end{equation}

Once we know the Hessian and the Macaulay polynomials associated to a nonsingular Hasse cubic $f_a=0$, we can determine if they coincide, or if they are projectively equivalent. 

The solutions of the equation $$\frac{4-a^3}{3a^2}=-\frac{2}{a}$$ correspond to the only cubics of the Hasse pencil having $H_a=M_a$ and $a \not = 0$: they are the ones corresponding to the values
$$
a=-2, a=1-\sqrt{3}, a=1+\sqrt{3}.
$$

 It is a well known fact that two planar cubics are projectively equivalent (or isomorphic) if and only if they have the same $j$-invariant.
Hence we can study the behavior of the Macaulay and the Hessian polynomials for a representative in each isomorphism class by means of the value of its $j$-invariant. The equation (\ref{fascio di Hasse}) allows us to make complete and explicit computations and, moreover, there is a formula for the $j$-invariant of a nonsingular cubic of that form (see \cite{RR}, Lemma 2.2):
\begin{eqnarray}\label{jinv}
j(f_a)=-\frac{a^3(a^3+8)^3}{(1-a^3)^3}.
\end{eqnarray}

\begin{proposition}
There are only four projective equivalence 
classes of smooth plane cubics such that the Macaulay and the Hessian polynomial are projectively equivalent. They are the ones corresponding to the values $0$, $64$,  $\omega$ and $\overline{\omega}$ of the $j$-invariant, where $\omega=352+96\sqrt{15}i$.
\end{proposition}
\begin{proof}
By the smoothness hypothesis $a^3 \not = 1$ in (\ref{fascio di Hasse}). Recalling that the Fermat cubic has $j$-invariant equal to $0$, we can furthermore suppose that in (\ref{fascio di Hasse}) $a \not = 0$, $a^3 \not = -8$. Using (\ref{ahessiano}),(\ref{amacaulay}),(\ref{jinv}) and imposing $j(H_a)=j(M_a)$, we get the following equation:
\begin{align}\label{concubi} 
(-8 - 20 a^3 + a^6)^2 \cdot  &
(262144 - 1212416 a^3 + 3248128 a^6 - 4353536 a^9 + \\ 
  & 3988672 a^{12} - 1649216 a^{15} + 248320 a^{18} - 656 a^{21} + a^{24})=0. \nonumber
\end{align} 
The solutions of (\ref{concubi}) represent  the cubics in the Hasse pencil having the Hessian and Macaulay polynomials projectively equivalent. 
Luckily we do not need to compute them explicitly, since we are just interested to know the value of their $j$-invariants. Given (\ref{jinv}), it is then enough to consider:
\begin{align}\label{senzacubi} 
(-8 - 20 \alpha +  \alpha^2)^2 \cdot  &
(262144 - 1212416  \alpha+ 3248128 \alpha^2 - 4353536  \alpha^3 + \\ 
  & 3988672  \alpha^{4} - 1649216  \alpha^{5} + 248320  \alpha^{6} - 656  \alpha^{7} +  \alpha^{8})=0, \nonumber
\end{align} 
where $\alpha=a^3$. Let $P$ be the first factor in (\ref{senzacubi}) and let $Q$ be the second. 
The equation $P=0$ can be easily solved in $\alpha$ and the associated $j$-invariants computed by means of (\ref{jinv}): only one value comes out, namely $j=64$.
The polynomial $Q$ has degree $8$ in $\alpha$ and $8$ distinct roots in $\mathbb{C}$. 
Again, we just need to compute $-\frac{\alpha(\alpha+8)^3}{(1-\alpha)^3}$ in every root of $Q$, i.e., we can consider $-\frac{\alpha(\alpha+8)^3}{(1-\alpha)^3} \in \frac{\mathbb{C}[\alpha]}{Q\mathbb{C}[\alpha]}$. Since $(1-\alpha)$ and $Q$ are coprime, then in  $\frac{\mathbb{C}[\alpha]}{Q\mathbb{C}[\alpha]}$ we have that 
\begin{align} \label{ridotto}
-\frac{\alpha(\alpha+8)^3}{(1-\alpha)^3} & = \frac{26476544}{19683} - \frac{82462208}{19683}\alpha + \frac{47739200}{6561}\alpha^2 - \frac{
 146525704}{19683}\alpha^3 + \\ & \frac{62168965}{19683}\alpha^4 - \frac{3144947}{6561}\alpha^5 + \frac{
 24926}{19683}\alpha^6 - \frac{38}{19683}\alpha^7. \nonumber
 \end{align}
Call $K$ the right-hand side of (\ref{ridotto}).  Since $K^2-704K+262144$ is divisible by $Q$, then the other two values of the $j$-invariant are the two complex solutions of $x^2-704x+262144=0$: $ j=\omega, j=\overline{\omega}$ where $\omega= 352+96\sqrt{15}i$.
 
In all remaining cases, that are infinitely many, $H_a$ and $M_a$ are not projectively equivalent. This concludes the proof.
\end{proof}

\begin{remark}
Given a nonsingular plane cubic curve $f$, notice that  $f$, $\textrm{Hess}(f)$ and $\textrm{Mac}(f)$ are all mutually projectively equivalent if and only if $j(f)=64$.
\end{remark}

\section{Plane quartics} \label{quartiche piane}
In this section we will make some remarks about the relationship between the Hessian and the Macaulay polynomials for smooth quartics $f$ in $\mathbb{P}^2$: without any pretence of completeness we just would like to bring forward some examples.

Since we do not have the $j$-invariant at our disposal - as opposed to the previous sections - we need other ways to find out if $\mathrm{Mac}(f)$ and $\mathrm{Hess}(f)$ are or are not projectively equivalent. First of all notice that, in general, a necessary condition for $\mathrm{Mac}(f)$ and $\mathrm{Hess}(f)$ to be projectively equivalent is that the two Hilbert functions associated to the apolar ring to $\mathrm{Mac}(f)$ and $\mathrm{Hess}(f)$, respectively, are exactly the same. We will soon see that this condition is not sufficient. 
Notice also that, since we are assuming $f$ smooth, then $J(f)$ is generated by a regular sequence, i.e., $R(f)$ is always a complete intersection. This is not at all the case for the apolar ring to the Hessian polynomial: even if $R/\mathrm{Hess}(f)^\perp$ is always Artinian Gorenstein (hence of dimension $0$), in most cases the minimal number of generators of $\mathrm{Hess}(f)^\perp$ exceeds $\dim(R)=3$. Therefore another necessary condition for $\mathrm{Mac}(f)$ and $\mathrm{Hess}(f)$ to be projectively equivalent is that $R/\mathrm{Hess}(f)^\perp$ is a complete intersection, i.e., $\mathrm{Hess}(f)^\perp$ must be generated by only $3$ polynomials. Unfortunately also in this case the condition is not sufficient.

Since a general quartic can be written as a sum of at most six powers of linear forms (see, for example, \cite{AH}), 
and - as already seen - we can just analyze one quartic for each class of projective equivalence, we have the following cases:

\begin{enumerate}
\item Fermat quartic: $f=x^4+y^4+z^4$. 
The Hessian polynomial coincides, up to a multiplicative constant, with the Macaulay polynomial: $x^2y^2z^2$ (this is a special case of Ex.\ \ref{fermat}).

\item Caporali quartics: $f=x^4+y^4+z^4+l(x,y,z)^4$, where $l(x,y,z)=ax+by+cz$ is a linear form, $a,b,c \in \C$ (see \cite[Ex.\ 6.11]{D}). 
In this case, depending on $a,b,c$, there are many examples for which the Hessian polynomial and the Macaulay polynomial cannot be projectively equivalent. Let us provide some explicit computations.

Case 1: Let $a=b=c=1$. 
In this case $\mathrm{Hess}^\perp$ must be generated by at least $9$ polynomials and therefore $R/\mathrm{Hess}^\perp$ is not a complete intersection, i.e., $\mathrm{Mac}(f)$ and $\mathrm{Hess}(f)$ cannot be projectively equivalent.

Case 2: Let $a=1,b=2,c=0$. In this case the two Hilbert functions $h_M$ and $h_H$ coincide and moreover $\mathrm{Hess}(f)^\perp$ is generated by three homogeneous polynomials of degree $3$: $ x^2y - (13/11)xy^2 - (5/11)y^3, x^3 - (12/11)xy^2 - (8/11)y^3, z^3$. Anyway it can be proved that $\mathrm{Mac}(f)$ and $\mathrm{Hess}(f)$ are not projectively equivalent: in fact $\mathrm{Mac}(f)= z^2(x^4-4x^3y-6x^2y^2+8x y^3-2y^4)$ and $\mathrm{Hess}(f)=z^2(x^4+4x^3y+(9/2)x^2y^2+xy^3+y^4)$, i.e., their null loci consist of a double line $\{z=0\}$ and four more lines passing through $[0,0,1]$ and intersecting $\{z=0\}$ in four distinct points. By means of cross-ratio and $j$-invariant it is then easy to see that these two configurations of four points on $\mathbb{P}^1$ are not projectively equivalent.

Case 3: Let $a=1,b=(-2)^{1/4},c=0$. In this case $\mathrm{Mac}(f)=z^2(x^4-(8/b) x^3y+(12/b^2) x^2y^2+4b xy^3-2y^4)$ and $\mathrm{Hess}(f)=z^2(b^2x^4+2b^3x^3y+2b xy^3+b^2y^4)$ are projectively equivalent: as in the previous case, we just need to see that the $j$-invariant of  the two quadruples of points is the same.
\item Clebsh quartics: $f=x^4+y^4+z^4+l_1(x,y,z)^4+l_2(x,y,z)^4$, where $l_1,l_2$ are linear forms (see \cite[Def.\ 6.12.1]{DK}). As before there are many examples for which $\mathrm{Hess}(f)$ and $\mathrm{Mac}(f)$ are not projectively equivalent. In some cases however (for example when $f=x^4+y^4+z^4+(x+y)^4+(x+2y)^4$) the apolar ring to $\mathrm{Hess}(f)$ is a complete intersection. 
\end{enumerate}

Given the fact, as we have just seen, that for a general quartic curve its Macaulay and Hessian polynomials are not projectively equivalent, it is quite surprising that for the Klein quartic $$f=x^3 y+y^3 z+ z^3 x,$$ its Hessian polynomial coincides, up to a multiplicative constant, with the Macaulay polynomial $\textrm{Mac}(f)=xy^5 + x^5z - 5x^2y^2z^2 + yz^5$. Notice that the same is true also for the Klein cubic in $\mathbb{P}^4$ $$ g= x^2y+y^2z+z^2w+w^2t+t^2x,$$ where $\textrm{Hess}(g)= 32 x^3z^2 - 32 xy z^3 + 32 y^3w^2 + 32 x^2w^3 - 32 yz w^3 -
32 xy^3t - 32 x^3 wt + 96 xy z wt + 32 z^3t^2 + 32 y^2t^3 - 32 zw t^3$. 

Notice that, unlike the Fermat case, the fact that $\textrm{Hess}(f)$ ($\textrm{Hess}(g)$) is projectively equivalent to $\textrm{Mac}(f)$ ($\textrm{Mac}(g)$, respectively) does not entirely depend on the particular ``symmetries'' of the arrangements of variables into the equation: for example for the cubic surface in $\mathbb{P}^3$ $$h=x^2y+y^2z+z^2w+w^2x,$$ we have that $\textrm{Hess}(h)$ is not projectively equivalent to $\textrm{Mac}(h)$, as it is clear, for example, computing the Hilbert functions of the two apolar rings.

\section{\cocoa}\label{cocoa}
The freely available \cocoa \ system, implemented by a team in Genoa (see \cite{Cocoa}), is well suited to perform calculations on polynomials. In particular it turned out to be very useful to experiment about the projective equivalence of the Hessian and the Macaulay polynomials associated to $f \in R_d$. This short program is an example on how we made calculations for Ex.\  \ref{esempio}:

\begin{itemize}
\item[] \texttt{Use R::=QQ[x,y,z];} defines the ring in which we want to work;\\
\item[] \texttt{F:=(x+y)\^{}3+y\^{}3+z\^{}3;} fixes the polynomial {$F$};\\
\item[] \texttt{N:=Jacobian([F]); H:=Jacobian(N[1]); G:=Det(H);} returns the Hessian polynomial of {$F$};\\
\item[]  \texttt{J:=Ideal(N[1]);} computes the Jacobian ideal of {$F$};\\
\item[] \texttt{Gort:=PerpIdealOfForm(G);} returns  the ideal of derivations that kill {$G$}, i.e. $G^\perp$;\\
\item[] \texttt{Hilbert(R/J);} computes the Hilbert function of {$R/J$};\\
\item[] \texttt{Hilbert(R/Gort);} computes the Hilbert function of {$R/G^\perp$};\\
\item[] \texttt{InverseSystem(J,3);} returns the Macaulay polynomial associated to {$F$} (that is, the Macaulay polynomial associated to {$R/J$}, which has degree $3$).\\
\end{itemize}

Also the function \texttt{DerivationAction(D,P);} is very useful: it returns the action of the derivation D on the polynomial P. For example, in Ex.\ \ref{esempio} \texttt{DerivationAction(216xyz+216y\^{}2z,3(x+y)\^{}2z);} returns \texttt{2592}. Finally \texttt{MinGens(I);} gives a minimal list of generators for the ideal \texttt{I}.

\section{Acknowledgments}
Both authors wish to warmly thank Prof.~ Edoardo Sernesi for having suggested this problem and for many helpful and fruitful discussions. The authors are also grateful to Prof.~ Ragni Piene and to the referee, whose suggestions improved the presentation of this paper.
\bibliographystyle{plain}
\bibliography{apolarity} 

\end{document}